\def\VersionDateTime{04/January/2017}

\def\yes{\if00}
\def\no{\if01}
\def\iftwelvept{\yes}

\def\ifusepdf{\no}
\def\ifpsfont{\yes}

\iftwelvept
\documentclass[leqno,12pt]{amsart}
\else
\documentclass[leqno]{amsart}
\fi
\usepackage{amscd}
\usepackage{amsmath}
\usepackage{amssymb}
\usepackage{amsfonts}
\usepackage{latexsym}
\usepackage[all]{xy}
\usepackage{graphicx}
\usepackage{multicol}
\usepackage{ascmac}
\usepackage{fancybox, framed}
\usepackage{verbatim}
\usepackage{mathrsfs}

\ifusepdf
\usepackage{hyperref}
\else\fi
\ifpsfont
\usepackage[T1]{fontenc}
\else\fi

\iftwelvept
\else\fi

\newtheorem{Theorem}{Theorem}[section]
\newtheorem{Proposition}[Theorem]{Proposition}
\newtheorem{Lemma}[Theorem]{Lemma}
\newtheorem{Corollary}[Theorem]{Corollary}

\theoremstyle{definition}
\newtheorem{Definition}[Theorem]{Definition}
\newtheorem{Remark}[Theorem]{Remark}
\newtheorem{Example}[Theorem]{Example}

\newcommand{\TT}{{\mathbb{T}}}
\newcommand{\GG}{{\mathbb{G}}}

\newcommand{\ZZ}{{\mathbb{Z}}}
\newcommand{\QQ}{{\mathbb{Q}}}
\newcommand{\RR}{{\mathbb{R}}}

\newcommand{\PP}{{\mathbb{P}}}
\newcommand{\OO}{{\mathcal{O}}}

\newcommand{\Lscr}{{\mathscr{L}}}
\newcommand{\Uscr}{{\mathscr{U}}}
\newcommand{\Sscr}{{\mathscr{S}}}
\newcommand{\Vscr}{{\mathscr{V}}}

\newcommand{\mfrak}{{\mathfrak{m}}}

\newcommand{\Irr}{\operatorname{\operatorname{Irr}}}

\newcommand{\ord}{\operatorname{ord}}

\newcommand{\rest}[2]{\left.{#1}\right\vert_{{#2}}}  

\newcommand{\Spec}{{\operatorname{Spec}}}

\newcommand{\trop}{\operatorname{trop}}

\newcommand{\an}{{\operatorname{an}}}

\usepackage{mathrsfs}
\newcommand{\Xscr}{{\mathscr{X}}}

\newcommand{\Proof}{{\sl Proof.}\quad}
\newcommand{\QED}{{\unskip\nobreak\hfil\penalty50\quad\null\nobreak\hfil
{$\Box$}\parfillskip0pt\finalhyphendemerits0\par\medskip}}

\begin{document}

\title[Tropicalization associated to adjoint linear systems]{Effective faithful tropicalizations associated to adjoint linear systems}
\author{Shu Kawaguchi}
\address{Department of Mathematical Sciences, 
Doshisha University, Kyotanabe, Kyoto 610-0394, Japan}
\email{kawaguch@math.kyoto-u.ac.jp}
\author{Kazuhiko Yamaki}
\address{Institute for Liberal Arts and Sciences, 
Kyoto University, Kyoto, 606-8501, Japan}
\email{yamaki.kazuhiko.6r@kyoto-u.ac.jp}
\date{\VersionDateTime}
\thanks{The first named author partially supported by KAKENHI 15K04817 and 25220701,  and the second named author partially supported by KAKENHI 26800012.}

\newcommand{\Proj}{\operatorname{\operatorname{Proj}}}
\newcommand{\Prin}{\operatorname{Prin}}
\newcommand{\Rat}{\operatorname{Rat}}
\newcommand{\zero}{\operatorname{div}}
\newcommand{\Sau}{\operatorname{Sau}}
\newcommand{\Func}{\operatorname{Func}}

\newcommand{\relint}{{\textup{relint}}}

\newenvironment{notation}[0]{%
  \begin{list}%
    {}%
    {\setlength{\itemindent}{0pt}
     \setlength{\labelwidth}{4\parindent}
     \setlength{\labelsep}{\parindent}
     \setlength{\leftmargin}{5\parindent}
     \setlength{\itemsep}{0pt}
     }%
   }%
  {\end{list}}

\newenvironment{parts}[0]{%
  \begin{list}{}%
    {\setlength{\itemindent}{0pt}
     \setlength{\labelwidth}{1.5\parindent}
     \setlength{\labelsep}{.5\parindent}
     \setlength{\leftmargin}{2\parindent}
     \setlength{\itemsep}{0pt}
     }%
   }%
  {\end{list}}
\newcommand{\Part}[1]{\item[\upshape#1]}

\begin{abstract}
Let $R$ be a complete discrete valuation ring of equi-characteristic zero with fractional field $K$. Let $X$ be a connected, smooth projective variety of dimension $d$ over $K$, and let $L$ be an ample line bundle over $X$. We assume that there exist a regular strictly semistable model $\Xscr$ of $X$ over $R$ and a relatively ample line bundle $\Lscr$ over $\Xscr$ with $\rest{\Lscr}{X} \cong L$. Let $S(\Xscr)$ be the skeleton associated to $\Xscr$ in the Berkovich analytification $X^{\an}$ of $X$. In this article, we study when $S(\Xscr)$ is faithfully tropicalized into tropical projective space by the adjoint linear system $|L^{\otimes m} \otimes \omega_X|$. Roughly speaking, our results show that, if $m$ is an integer such that the adjoint bundle is basepoint free, then the adjoint linear system admits a faithful tropicalization of  $S(\Xscr)$. 
\end{abstract}

\maketitle

%
\setcounter{Theorem}{0}
\section{Introduction}
\label{sec:intro}
Let $X$ be a connected, smooth projective variety defined over a field $K$, let $L$ be an ample line bundle over $X$, and let $\omega_X$ be the canonical line bundle over $X$. If $X$ is one dimensional and $\deg(L) \geq 2g(X) + 1$, where $g(X)$ is the genus of $X$, then global sections of $L$ give an embedding of $X$ into projective space. 
In higher dimension, adjoint bundles $L^{\otimes m} \otimes \omega_X$ can be viewed as the generalization of line bundles of large degree on curves, 
and there have been many important studies on when global sections of $L^{\otimes m} \otimes \omega_X$ give an embedding (or injection) of $X$ into projective space (see, for example, \cite[II, 10.4]{La}). 

We consider an analogue for Berkovich skeleta in non-Archimedean geometry, and we study when global sections of line bundles give a faithful tropicalization of a Berkovich skeleton into tropical projective space. In \cite{KY}, we have studied a one-dimensional case. In this paper,
with an assumption on the existence of a model, we consider a faithful tropicalization by global sections of adjoint bundles in arbitrary dimension.  

In the following, we assume that $K$ is a complete discrete valuation field. 
We denote by $R$ the ring of integers and by $k$ the residue field. We assume that $R$ has 
equi-characteristic zero and that $X$ admits regular strictly semistable model $\Xscr$ of $X$ over $R$. (If we replace $R$ by a suitable a finite extension, then there always exists such $\Xscr$, see \S\ref{subsec:semistable} for details.) Let $X^\an$ denotes the analytification of $X$ in the sense of Berkovich. Let $S(\Xscr) \subset X^{\an}$ be the Berkovich skeleton associated to $\Xscr$, which carries a piecewise integral rational-affine structure and is homeomorphic to the dual intersection complex of the special fiber $\Xscr_s$  (see \S\ref{subsec:skeleton} for details). 

We say that $S(\Xscr)$ has a {\em faithful tropicalization associated to the linear system $|L|$} if there exist nonzero global sections $s_0, \ldots, s_n \in H^0(X, L)$ with the following properties: 
\begin{equation}
\label{eqn:trop:map:L}
  \varphi: X^\an \longrightarrow \TT\PP^n, \quad p = (p, |\cdot|) \mapsto \left(-\log|s_0(p)|: \cdots: -\log|s_n(p)|\right) 
\end{equation}
is a morphism into tropical projective space $\TT\PP^n$ and  the restriction of 
$\varphi$ to $S(\Xscr)$ is a homeomorphism onto its image preserving the piecewise integral rational-affine structures. (For the tropical projective space and faithful tropicalizations, see \S\ref{subsec:tropical} and \S\ref{subsec:unimodular:faithful} for details.)

\medskip
For a positive integer $d$, we set 
\begin{equation}
\label{eqn:def:phi:d}
\phi(d) 
:= \min \left\{
m_0 \in \ZZ \;\left|\;  
\begin{aligned}
& \text{For any smooth projective variety $Z$ with $\dim Z \leq d$ over} \\
& \text{a field of characteristic zero, and for any ample line bundle}\\
& \text{$N$ over $Z$, $N^{\otimes m} \otimes \omega_Z$ is basepoint free  for any $m \geq m_0$.}
\end{aligned}
\right.
\right\}, 
\end{equation}
where $\omega_Z$ denotes the canonical line bundle of $Z$. 
By Angehrn and Siu \cite{AS}, 
one has $\phi(d) \leq d(d+1)/2+1$, and by a subsequent work of Heier \cite{Heier},  one has $\phi(d) \leq O(d^{4/3})$. Fujita's conjecture concerns the optimal $\phi(d)$, which asks if $\phi(d) = d+1$ holds true, and is known to be true in dimension $d = 2, 3, 4$ (see Reider~\cite{Re}, Ein--Lazarsfeld~\cite{EL}, and Kawamata~\cite{Ka}). 

\smallskip
Let us state our main theorem. 

\begin{Theorem}
\label{thm:main}
Let $R$ be a complete discrete valuation ring of equi-characteristic zero with fractional field $K$. Let $X$ be a connected, smooth projective variety of dimension $d$ over $K$, and let $L$ be an ample line bundle over $X$. Assume that $X$ has a regular strictly semistable model  $\Xscr$ over $R$ such that $L$ extends to a relatively ample line bundle $\Lscr$ over $\Xscr$. 
Let $S(\Xscr)$ denote the associated skeleton in $X^{\an}$. Then $S(\Xscr)$ has a faithful tropicalization associated to $\left| L^{\otimes m} \otimes \omega_X \right|$ for any $m \geq \phi(d)$. More strongly, if $\ell$ denotes the number of irreducible components of $\Xscr_s$, 
then there exist $(\ell+d+1)$ nonzero global sections of $L^{\otimes m} \otimes \omega_X$ such that the associated morphism $\varphi: X^{\an} \to \TT\PP^{\ell+d}$ is a faithful tropicalization of $S(\Xscr)$. 
\end{Theorem}

\begin{Corollary}
\label{cor:faithful:canonical}
Let $R$ and $K$ be as in Theorem~\ref{thm:main}. 
Let $X$ be a connected, smooth projective variety of dimension $d$ over $K$, and 
assume that $X$ has a regular strictly semistable model $\Xscr$ over $R$. 
Let $S(\Xscr)$ denote the the associated skeleton in $X^{\an}$. 
\begin{enumerate}
\item
Let $L$ be an ample line bundle over $X$, and we assume that $L$ extends to a relatively ample line bundle over $\Xscr$. Then, if $\omega_X$ is trivial \textup{(}e.g. if $X$ is an abelian variety or a Calabi-Yau manifold\textup{)} and if $d \leq 4$, then  $S(\Xscr)$ has a faithful tropicalization associated to $\left| L^{\otimes m} \right|$ for any $m \geq d + 1$. 
\item
If $\omega_{\Xscr/R}$ is relatively ample and if $d \leq 4$, then  $S(\Xscr)$ has a faithful tropicalization associated to $\left| \omega_X^{\otimes m} \right|$ for any $m \geq d + 2$. 
\end{enumerate}
\end{Corollary}

Faithful tropicalizations have been studied by many authors.
Following the one-dimensional case due to
Katz--Markwig--Markwig \cite{KMM, KMM2}, Baker--Payne--Rabinoff \cite{BPR}, 
and Baker--Rabinoff \cite{BR},
Gubler--Rabinoff--Werner \cite{GRW} have shown that for a connected smooth projective variety $X$ of an arbitrary dimension
and for any skeleton $\Gamma$ of $X^{\an}$, there exist nonzero rational functions $f_1, \ldots, f_n$ of $X$ such that
\[
  \psi: X^{\an} \dasharrow \RR^n, \quad p = (p, |\cdot|)\mapsto (-\log|f_1(p)|, \ldots, -\log|f_n(p)|)
\] 
gives a faithful tropicalization of $\Gamma$. A new aspect here is to give a condition on when $f_1, \ldots, f_n$ are taken to be global sections of the adjoint line bundle. (See also \cite{KY} for related results when $\dim X = 1$). 

While the assumption on the existence of a model in Theorem~\ref{thm:main} may look strong, we remark that, 
if we replace $K$ with a finite extension field of $K$ and $L$ with a power of $L$,
then we can always achieve this assumption, 
so that we can apply the theorem; see Remark~\ref{rmk:learned:from:Gubler}.

Musta\k{t}\u{a}--Nicaise \cite{MN} and Nicase--Xu \cite{NX1} have introduced the essential skeleton $S(X)$ which depends only on $X$ and not on a particular model $\Xscr$ of $X$
by using birational geometry. 
If $X$ has a regular strictly semistable model of $\Xscr$, 
then $S(X)$ is contained in $S(\Xscr)$, so that Theorem~\ref{thm:main} also gives a faithful tropicalization of $S(X)$. 

It will be interesting to study a faithful tropicalization of the (more general) skeleton $S(\Xscr, H)$ associated to a strictly semistable pair (cf. \cite{GRW}) or the skeleton associated to an sncd-model of $X$ (cf. \cite{BFJ}, \cite{MN}). 
It may be possible to work over an  algebraically closed field which is complete with respect to a non-trivial non-Archimedean value, but working under this setting should be technically more involved (cf. \cite{KY}), and we do not pursue it in this paper. 

Let us explain our ideas of the proof of Theorem~\ref{thm:main}, together with the structure of this paper. In Section~\ref{sec:prelim}, we briefly recall some known facts on Berkovich spaces and tropical geometry. 
In Section~\ref{sec:lemmas}, we prove 
some vanishing of cohomologies and basepoint-freeness on strictly normal
crossing varieties, which are technical keys.
Then in Section~\ref{sec:proof:unimodular}, we start the proof of Theorem~\ref{thm:main}. 
We denote by $\Xscr_s  = \bigcup_{1 \leq i \leq \ell}\Xscr_{s, i}$ the irreducible decomposition of the special fiber $\Xscr_s$ of $\Xscr$. Then each $\Xscr_{s, i}$ corresponds to a point, called a Shilov point
 (also called a divisorial point) 
and denoted by $[\Xscr_{s, i}]$,  in $S(\Xscr)$. We need, for example, to separate $[\Xscr_{s, i}]$ and $[\Xscr_{s, i^\prime}]$ ($i \neq i^\prime$) by global sections of $L^{\otimes m} \otimes \omega_X$. To this end, we first construct suitable 
global sections $\eta$ of $\rest{\Lscr^{\otimes m}\otimes \omega_{\Xscr/R}}{\Xscr_{s, i}}(-W)$ for suitable effective divisors $W$ of $\Xscr_{s, i}$ when $m \geq \phi(d)$. Next, 
we extend those $\eta$'s to global sections of $\rest{\Lscr^{\otimes m}\otimes \omega_{\Xscr/R}}{\Xscr_{s, i}}$, and 
using a Kodaira-type vanishing of cohomologies, we extend them to global sections $\widetilde{\eta}$ of 
$\Lscr^{\otimes m}\otimes \omega_{\Xscr/R}$, where $\omega_{\Xscr/R}$ is the relative dualizing sheaf. 
Then we restrict to the $\widetilde{\eta}$'s to the generic fiber $X$ to obtain global sections $s$ of $L^{\otimes m} \otimes \omega_X$, and we show that 
those $s$'s give a local homeomorphism preserving the integral structures between $S(\Xscr)$ and its image in tropical projective space, which is called a unimodular tropicalization of $S(\Xscr)$.  
Actually, since we only consider the skeleton $S(\Xscr)$ associated to a regular strictly semistable model $\Xscr$, we show that $(\ell+d
+1)$ nonzero global sections $s_0, s_1, \ldots, s_{\ell+d
}$, where $\ell$ is the number of irreducible components of $\Xscr_s$, are enough to give a unimodular tropicalization. 
In Section~\ref{sec:proof:faithful}, we show somewhat fortunately that this $\varphi$ is a faithful tropicalization of $S(\Xscr)$.  

\smallskip
{\sl Acknowledgments.}\quad 
We deeply thank Professor Osamu Fujino for helpful discussions on Kodaira-type vanishing theorems, and 
Professor Walter Gubler for helpful discussions on Remark~\ref{rmk:learned:from:Gubler}. 

\subsection*{Notation and conventions}
Throughout this paper, let $R$ be a complete discrete valuation ring of equi-characteristic zero with fractional field $K$ and residue field $k$. Let $\varpi$ denote a uniformizer of $R$, and let $\mfrak$ denote the maximal ideal of $R$.  
Let $v_K: K \to \ZZ \cup \{\infty\}$ denote the additive discrete valuation of $K$ with $v_K ( \varpi) = 1$, and let 
 $|\cdot|_K := \exp(-v_K(\cdot))$ denote the multiplicative value of $K$. 
 
By a {\em variety}, we mean a reduced separated scheme 
of finite type over a field, which is allowed to be reducible.

\section{Preliminaries}
\label{sec:prelim}
In \S\S\ref{subsec:semistable}--\ref{subsec:tropical}, we collect some known facts on Berkovich spaces and tropical geometry. Our basic references are \cite{BPR}, \cite{BFJ},  \cite[\S5]{Gu4}, \cite{GRW}, \cite{MS}, \cite{MN} and \cite{N1}. Similar to  \cite{BFJ} and \cite{MN}, we use the language of schemes (rather 
than formal schemes), which might be more familiar to readers in birational algebraic geometry. 
Then, in \S\ref{subsec:unimodular:faithful}, following \cite{BPR} and \cite{GRW}, 
we define unimodular and faithful tropicalizations associated to linear systems. 

\subsection{Strictly semistable model}
\label{subsec:semistable}
Let $X$ be a connected, smooth projective variety defined over $K$ of dimension $d$. 
By a {\em strictly semistable model} $\Xscr$ of $X$, 
we mean a scheme that is proper and flat over $R$  endowed with an isomorphism 
$\Xscr\times_{\Spec(R)}\Spec(K) \cong X$ 
such that $\Xscr$ is covered by Zariski open subsets $\Uscr$ that admit an {\'e}tale morphism
\begin{equation}
\label{eqn:etale:mor}
\psi: \Uscr \to \Sscr := \Spec\, R[x_0, \ldots, x_d]/(x_0 x_1\cdots x_r - \varpi), 
\end{equation}
where $r$ is an integer with $0 \leq r \leq d$ 
depending on $\mathscr{U}$ (cf. \cite[(2.16)]{deJong}, \cite[Definition~3.1]{GRW}).
Remark that any irreducible component
of the special fiber of a strictly semistable model
is irreducible.
If  a strictly semistable model $\Xscr$ is regular, 
we call $\Xscr$  a {\em regular} strictly semistable model of $X$.

We note that $X$ may not admit a regular strictly semistable model over $R$, but after 
making a suitable finite extension, one can always find a regular strictly semistable model of $X$.
Indeed, by Hironaka's resolution of singularities (recall that we assume that $R$ has equi-characteristic zero), one can find a proper regular scheme 
$\Xscr$ over $R$ with generic fiber $X$ for which the reduced closed fiber $\Xscr_{s,{\rm red}}$ 
is a divisor with normal crossings in $\Xscr$. Then the semistable theorem by 
Kempf--Knudsen--Mumford--Saint-Donat  \cite[p. 198]{KKMS} asserts that, 
after making a finite extension $K^\prime/K$ with discrete valuation ring $R^\prime$, 
$X \times_{\Spec(K)} \Spec(K^\prime)$ admits 
a regular strictly semistable model over $R^\prime$.

\subsection{Berkovich analytic space}
\label{subsec:Berko}
Let $X$ be an irreducible variety defined over $K$. 
Then one has 
$X^{\an}$,
the {\em analytification} of $X$ in the sense of Berkovich 
(cf. \cite{Be1}, \cite{Be2}, \cite{Be3}, \cite{Be4}). 
Although $X^{\an}$ has an analytic structure, 
we recall only how it is described  as a topological space.
For a scheme point $p$ in $X$, let $\kappa(p)$ denote the residue field. 
As a set, $X^{\an}$ is given by 
\[
  X^{\an} := 
  \left\{
  (p, |\cdot|) \mid \text{$p \in X$ and $|\cdot|$ is an absolute value of $\kappa(p)$ extending $|\cdot|_K$}
  \right\}. 
\]
The underlying topology on $X^{\an}$ is the weakest topology such that 
$\iota: X^\an \to X, (p, |\cdot|) \mapsto p$ is continuous and that for any Zariski open set $U$ of $X$ and for any regular function $g \in \OO_X(U)$, the map  $\iota^{-1}(U) \to \RR, (p, |\cdot|) \mapsto |g(p)|$ is continuous. 

\begin{Example}[classical point]
\label{eg:classical}
If $p$ is a closed point of $X$, then $\kappa(p)$ is a finite extension 
field of the completely valued field $K$, and 
hence there exists a unique 
absolute value $|\cdot|_{\kappa(p)}$ that extends to $|\cdot|_K$. Through the assignment $p \mapsto (p, |\cdot|_{\kappa(p)})$, any closed point of $X$ is regarded as an element of $X^{\an}$. 
\end{Example}

\begin{Example}[Shilov point]
Let $\Xscr \to \Spec(R)$ be a proper 
and flat morphism with generic fiber $X$. 
Assume that $\Xscr$ is normal.
Let $E$ be an irreducible component of the special fiber, 
and let $\xi$ denote the generic point of $E$. 
Then since $\xi$ is a normal point of codimension $1$,
$\OO_{\Xscr , \xi}$ is a discrete valuation ring.
Further,
$\OO_{\Xscr , \xi}$
contains $R$
and
the fraction field of $\OO_{\Xscr , \xi}$ equals the function field
$K(X)$ of $X$.
It follows that there exists a unique value 
$|\cdot|_\xi$ on $K(X)$ 
that is equivalent to the value associated to
the discrete valuation ring
$\OO_{\Xscr , \xi}$ and that extends $|\cdot|_K$.
Thus, if $\eta$ denotes the generic point of $X$, then 
$(\eta, |\cdot|_\xi) \in X^{\an}$. The point  $(\eta, |\cdot|_\xi)$, simply denoted by $[E]$, 
is called the {\em Shilov point} associated to $(\Xscr, E)$. 
Shilov points are also called divisorial points
(cf. \cite{BFJ}, \cite{N1}, for example).

In particular, if $X$ is a connected, smooth projective variety and admits a strictly semistable model $\Xscr$ over $R$, then each irreducible component of the special fiber $\Xscr_s$ gives a Shilov point in $X^{\an}$. 
\end{Example}

\subsection{Skeleton}
\label{subsec:skeleton}
Let $X$ be a connected, smooth projective variety defined over $K$ of dimension $d$, and 
assume that $X$ admits a regular strictly semistable model $\Xscr$ over $R$. 
Then one has the skeleton $S(\Xscr) \subset X^{\an}$ associated to $\Xscr$, 
as a special case of \cite{Be3} (see also \cite{GRW}, \cite{KS}, \cite{N0}, \cite{T}). 
Here we mainly follow 
\cite[\S3]{BFJ} and \cite[\S3]{MN} to describe $S(\Xscr)$. 

\subsubsection*{Stratification}
Let $\Xscr_s$ denote the special fiber of $\Xscr$. 
We write ${\Xscr}_s  = \bigcup_{1 \leq i \leq \ell
}\Xscr_{s, i}$ for the irreducible decomposition. 
By the definition of a regular strictly semistable model, 
${\Xscr}_s$ is a reduced divisor of $\Xscr$, and each $\Xscr_{s, i}$ is a connected, smooth projective variety over $k$. 

We define the stratification of ${\Xscr}_s$ as follows.
We put ${\Xscr}_s^{(0)} := {\Xscr}_s$. For each $\alpha \in \ZZ_{> 0}$, let 
${\Xscr}_s^{(\alpha)}$ be the complement of the set of normal points in ${\Xscr}_s^{(\alpha-1)}$. Thus we obtain a chain of closed subsets: 
\[
{\Xscr}_s = {\Xscr}_s^{(0)}  \supsetneq 
{\Xscr}_s^{(1)}  \supsetneq 
\cdots \supsetneq 
{\Xscr}_s^{(t)}  \supsetneq 
{\Xscr}_s^{(t+1)}  = \emptyset, 
\]
where $t \leq d$. An irreducible components of ${\Xscr}_s^{(\alpha)}\setminus  {\Xscr}_s^{(\alpha+1)}$ ($0 \leq \alpha \leq t$) is called a {\em stratum} of ${\Xscr}_s$. A stratum $S$ of ${\Xscr}_s$ is {\em minimal} if there does not exist a stratum of ${\Xscr}_s$ that is strictly contained in 
the Zariski closure of $S$.

\subsubsection*{Canonical simplex associated to a stratum} 
Let $S$ be a stratum of $\Xscr$. 
Then there exist $r \geq 1$ and 
a subset $J = \{j_1, \ldots, j_r\}$ of $\{1, \ldots, \ell\}$ such that $S$ is a connected component of $\left(\bigcap_{j \in J} \Xscr_{s, j}\right) \setminus  {\Xscr}_s^{(r)}$. 

We define the standard $(r-1)$-simplex $\Delta^{r-1}$ by 
\begin{equation}
\label{eqn:standard:simplex}
\Delta^{r-1} := 
\{
\mathbf{u} := (u_1, \ldots, u_r) \in \RR^{r} \mid 
u_1 \geq 0, \ldots, u_r \geq 0, \, 
u_1 +  \cdots + u_r = 1
\}. 
\end{equation}
Let $\relint(\Delta^{r-1}) := 
\{
\mathbf{u} := (u_1, \ldots, u_r) \in \RR^{r} \mid 
u_1 > 0, \ldots, u_r > 0, \, 
u_1 +  \cdots + u_r = 1
\}
$ denote the relative interior of $\Delta^{r-1}$. 

For any $\mathbf{u} \in \Delta^{r-1}$, we are going to define an absolute value $|\cdot|_{\mathbf{u}, S}$ on the function field $K(X)$ extending $|\cdot|_K$. 
Let $\xi$ be the generic point of $S$. For each $j \in J$, we choose a local equation $T_j = 0$ of $\Xscr_{s, j}$ in $\Xscr$ at $\xi$. 
We take any $f \in \OO_{\Xscr, \xi}$. 
Let $\widehat{\OO}_{\Xscr, \xi}$ be the completion of $\OO_{\Xscr, \xi}$. 
Since $R$ is of equi-characteristic zero, Cohen's structure theorem asserts that 
$\widehat{\OO}_{\Xscr, \xi}$ contains a field isomorphic to 
the residue field $\kappa(\xi)$ and 
$\widehat{\OO}_{\Xscr, \xi}$ is isomorphic to the power 
series ring $\kappa(\xi)[\![T_{j_1}, \ldots, T_{j_r}]\!]$. 
Then 
$f$ is written
in $\widehat{\OO}_{\Xscr, \xi} \cong \kappa(\xi)[\![T_{j_1}, \ldots, T_{j_r}]\!]$ as 
\[
  f = \sum_{\mathbf{m} = (m_1, \ldots, m_r) \in {\ZZ_{\geq 0}}^r} a_{\mathbf{m}} T_{j_1}^{m_1} \cdots T_{j_r}^{m_r} 
  \quad\quad\quad
  (a_{\mathbf{m}} \in \kappa(\xi)), 
\]
and
\begin{equation}
\label{def:mathu:S}
  |f|_{\mathbf{u}, S} := \max \left\{ \exp(-u_1 m_1 - \cdots - u_r m_r) \mid \mathbf{m}  \in  {\ZZ_{\geq 0}}^r, 
  a_{\mathbf{m}}  \neq 0 \right\}
\end{equation}
gives a well-defined absolute value on $ \OO_{\Xscr, \xi}$ (see \cite[Proposition~3.1.4]{MN}). This absolute value $|\cdot|_{\mathbf{u}, S}$ extends to an absolute value $|\cdot|_{\mathbf{u}, S}$ on $K(X)$. 
Further, since $\varpi = \lambda T_{j_1} \cdots T_{j_r}$
for some unit $\lambda$ in $ \OO_{\Xscr, \xi}$, one has $|\varpi|_{\mathbf{u}, S} = \exp(-1)$,
which shows that
$|\cdot|_{\mathbf{u}, S}$ agrees with $|\cdot|_K$ for elements in $K$. 
Thus for each $\mathbf{u}$, we have a point $( \eta, |\cdot|_{\mathbf{u}, S} )
\in X^{\an}$, where $\eta$ is the generic point of $X$.

With the above notation,
we set 
\begin{equation}
\label{eqn:coord:DS}
  \Delta_S := \{(\eta, |\cdot|_{\mathbf{u}, S}) \in X^{\an} \mid \mathbf{u} \in \Delta^{r-1}\}
\end{equation}
and call it the {\em canonical simplex} associated to $S$. The assignment 
\begin{equation}
\label{eqn:intrinzic}
  \Delta^{r-1} \to \Delta_S, \quad \mathbf{u} \mapsto (\eta, |\cdot|_{\mathbf{u}, S}) 
\end{equation}
is a homeomorphism, where $\Delta^{r-1}$ is endowed 
with the Euclidean topology (see \cite[Proposition~3.1.4]{MN}). 
Via this homeomorphism,
we endow $\Delta_S$ with a simplex structure.
We denote $\relint(\Delta_S)$
the relative interior of the simplex $\Delta_S$.
Remark that
$
 \relint(\Delta_S) \cong
\relint(\Delta^{r-1})$ under the homeomorphism. 

Note that by the construction, 
$|\cdot|_{\mathbf{u}, S}$ is independent of the choice of 
local equations $T_j = 0$ for $\Xscr_{s, j}$. Thus, up to permutations of $J = \{j_1, \ldots, j_r\}$, 
the homeomorphism between $\Delta^{r-1}$ and $\Delta_S$ in 
\eqref{eqn:intrinzic} is intrinsic (cf. \cite[p.~169]{GRW}). 

Let $S^\prime$ be another stratum of which Zariski closure 
$\overline{S^\prime}$
contains $S$. 
Then
$\Delta_{S^\prime}$ is a face of $\Delta_S$, as we now explain.
Since $\overline{S^\prime} \supsetneq S$,
there exists a nonempty proper subset $A$ of 
$\{1, \ldots, r\}$ such that 
$S^\prime$
is defined in $\mathscr{X}$
at its generic point by 
$T_{j_\alpha} = 0$ for all $\alpha \in A$.
We set $r^\prime := |A|$. We regard the standard $(r^\prime-1)$-simplex 
$\Delta^{r^\prime-1}$ as the subset of $\Delta^{r-1}$ by 
\[
  \Delta^{r^\prime-1} = \{\mathbf{u} = (u_1, \ldots, u_r) \in \Delta^{r-1} \mid \text{$u_{\beta} = 0$ for any 
  $\beta \not\in A$}\}. 
\]
Then, by the definition of $|\cdot|_{\mathbf{u}, S}$, we see that 
$|\cdot|_{\mathbf{u}, S} = |\cdot|_{\mathbf{u}, S^\prime}$ for any $\mathbf{u} \in \Delta^{r^\prime-1}$. 
It follows that 
the homeomorphism 
$\Delta^{r^\prime-1} \to \Delta_{S^\prime}$ for $S^\prime$ coincides with the restriction to $\Delta^{r^\prime-1}$ of the homeomorphism $\Delta^{r-1} \to \Delta_S$ for $S$. 

\subsubsection*{Skeleton $S(\Xscr)$}
The {\em skeleton} $S(\Xscr) \subset X^{\an}$ associated to $\Xscr$ is defined by  
\[
  S(\Xscr) := \bigcup_{S} \Delta_S,
\]
where $S$ runs through all the strata of $\Xscr_s$. 

The description of the canonical simplex associated to a stratum tells us that 
$S(\Xscr)$ is homeomorphic to the {\em dual intersection complex} of $\Xscr_s$. The dual intersection complex of $\Xscr_s$ is a simplicial complex whose simplices 
correspond bijectively to the set of strata of $\Xscr_s$. To be precise, simplices of dimension $\alpha$ correspond 
bijectively to the irreducible components of ${\Xscr}_s^{(\alpha)}\setminus  {\Xscr}_s^{(\alpha+1)}$. 

As explained above, if $S$ and $S^\prime$ are strata of $\Xscr_s$, then the simplex corresponding to $S^\prime$ is a face of the simplex corresponding to $S$ if and only if the Zariski closure of $S^\prime$ contains $S$. In particular, vertices (i.e., $0$-dimensional simplices) of the dual intersection complex correspond bijectively to the irreducible components 
of $\Xscr_s$. The assignment of the simplex corresponding to a stratum $S$ to $\Delta_S$ in $X^{\an}$ gives 
the homeomorphism between the dual intersection complex of $\Xscr_s$ and 
the skeleton $S(\Xscr)$. See \cite{KS}, \cite[\S3]{BFJ}, \cite[\S3]{MN} for details. 

We note that 
$
  S(\Xscr) = \bigcup_{S} \Delta_S,
$
where $S$ runs through all the {\em minimal} strata of $\Xscr_s$, and that 
$
  S(\Xscr) = \coprod_{S} \relint(\Delta_S),
$
where $S$ runs through all the  strata  of $\Xscr_s$. 

\begin{Remark}
\label{rmk:learned:from:Gubler}
Suppose that we are given a skeleton $S ( \Xscr )$,
where
$\Xscr$ is a
strictly semistable model of $X$.
Then modulo finite extensions of $K$,
we can faithfully tropicalize $S ( \Xscr )$
by using Theorem~\ref{thm:main}.
Indeed, 
if we replace $K$ by a suitable finite extension $K^{\prime}$ with discrete valuation ring $R^{\prime}$ and $L$ by a suitable power $L^{\otimes a}$, then 
$X^{\prime} := X \times_{\Spec(K)} \Spec(K^{\prime})$ admits a regular strictly semistable model  $\Xscr^{\prime}$ over $R^{\prime}$ such that $S(\Xscr) \subseteq S(\Xscr^{\prime})$ (in the Berkovich analytification $X^{\prime \an}$) and such that $L^{\prime} := L^{\otimes a} \otimes_K K^{\prime}$ extends to a relatively ample line bundle over $\Xscr^{\prime}$. 
(Here, we use the same argument as in \cite[Corollary~1.5]{BFJ}.)
This means that
modulo finite extensions of $K$ and powers of $L$, 
there always exists a model $(\Xscr^{\prime}, \Lscr^{\prime})$ 
such that $\Xscr^{\prime}$ is strictly semistable,
$\Lscr^{\prime}$ is relatively ample
and such that
$S(\Xscr) \subseteq S(\Xscr^{\prime})$.
Thus Theorem~\ref{thm:main}
gives a faithful tropicalization of
$S (\Xscr^{\prime})$,
which means that
modulo finite extensions of $K$,
the theorem
gives a faithful tropicalization of
$S (\Xscr)$.
\end{Remark}

\subsection{Tropical projective space and tropicalization map}
\label{subsec:tropical}
In this subsection, we recall tropical projective space and tropicalization map from a Berkovich analytic space to 
tropical projective space. 

\subsubsection*{Integral rational-affine polyhedra}
Let $N$ be a free abelian group of rank $n$,
 and set $N_\RR := N \otimes_\ZZ \RR$. 
By fixing a free basis of $N$, we identify $N$ with $\ZZ^n$ and $N_\RR$ with Euclidean space $\RR^n$. 
Let $\langle\;,\;\rangle$ denote the standard inner product of $\RR^n$. 
An {\em integral rational-affine polyhedron} in $N_\RR = \RR^n$ 
is a subset written as 
$\Delta := \{x \in \RR^n \mid A x \leq b\}$ 
for some $A \in M_{r, n}(\ZZ)$ and $b \in \QQ^r$. 
The {\em face} of an integral rational-affine polyhedron $\Delta$ is a subset 
$\{x \in \Delta \mid \text{$\langle w, x\rangle \leq  \langle w, y\rangle$ for all $y \in \Delta$}\}$ 
for some $w \in \ZZ^n$. The relative interior of $\Delta$ is denoted by $\relint(\Delta)$. 
An {\em integral rational-affine polyhedral complex} $\Sigma$ is a finite union of 
integral rational-affine polyhedra such that any face of $\Delta \in \Sigma$ belongs to $\Sigma$ and such that 
any nonempty intersection of two polyhedra $\Delta, \Delta^\prime \in \Sigma$ is
a face of both $\Delta$ and $\Delta^\prime$. 
An {\em integral rational-affine map} from $\RR^n$ to $\RR^{m}$ is a map of the form 
$
  x \mapsto C x + d
$ for some $C \in M_{m, n}(\ZZ)$ and $d \in \QQ^m$. 

Let $\Delta \subset \RR^n$ be an integral rational-affine polyhedron. A 
map $F: \Delta \to \RR^m$ is an {\em integral rational-affine map} if $F$ is the restriction of an 
integral rational-affine map from $\RR^n$ to $\RR^{m}$. 
Let $H$ be the linear span of $\Delta-x$ for any $x \in \Delta$. 
Then $H \cap \ZZ^n$ is a free abelian group of rank $\dim \Delta$. 
If $F: \Delta \to \RR^m$ is an integral rational-affine map, then 
$\Delta^\prime := F(\Delta)$ is  an integral rational-affine polyhedron in $\RR^m$. 
Similarly, if $H^\prime$ is the linear span of $\Delta^\prime-x^\prime$ for any $x^\prime 
\in \Delta^\prime$ then $H^\prime \cap \ZZ^m$ is a free abelian group 
of rank $\dim \Delta^\prime$. We say that $F$ is a {\em unimodular} integral rational-affine map 
if $F$ is injective and $H^\prime \cap \ZZ^m = F(H \cap \ZZ^n)$ (see \cite[\S2.2]{GRW}). 

\subsubsection*{Tropical projective space}
We set  $\TT := \RR\cup\{\infty\}$. The $n$-dimensional tropical projective space is defined to be 
\[
\TT\PP^n := 
(\TT^{n+1}\setminus\{(\infty, \ldots, \infty)\})/\sim, 
\]
where by definition
$x := (x_0, \ldots, x_n), y := (y_0, \ldots, y_n) \in \TT^{n+1}\setminus\{(\infty, \ldots, \infty)\}$ satisfy $x \sim y$ if 
and only if there exists $c \in \RR$ such that $y_i = x_i + c$ for all $i = 0, \ldots, n$ (see \cite{MZ}). The equivalence class of $x$ in $\TT\PP^n$ is written as 
$(x_0: \cdots: x_n)$. Tropical projective space is equipped with $(n+1)$ charts 
$U_i := \{x = (x_0: \cdots: x_n) \in \TT\PP^n\mid x_i \neq \infty\}$. 

\subsubsection*{Tropicalization map}
Let $X_0, \ldots, X_n$ be the homogeneous coordinates of $\PP^n$. Let $\PP^{n, \an}$ denote the Berkovich analytification of $\PP^n$. 
Then we have the tropicalization map 
\begin{equation}
\label{eqn:trop:map}
\trop: \PP^{n, \an} \to \TT\PP^n, \qquad p = (p, |\cdot|)  \mapsto (-\log|X_0(p)|
: \cdots : -\log|X_n(p)|). 
\end{equation}

Since we fix the homogeneous coordinates, the multiplicative group $\GG_m^n$ is embedded in $\PP^{n}$, where $\GG_m^n = \{(1: x_1: \cdots: x_n) \mid x_1 \neq 0, \ldots, x_n \neq 0 \}$. Similarly, the Euclidean space $\RR^n$ is embedded in $\TT\PP^n$, where 
$\RR^n = \{(0: x_1: \cdots: x_n) \mid x_1, \ldots, x_n \in \RR\}$. 
Then the restriction of the  tropicalization map to $\GG_m^n$ gives 
$\rest{\trop}{\GG_m^n}: \GG_m^n \to \RR^n$,
and $\trop^{-1} ( \RR^{n} ) = \GG_m^n$.

Tropical geometry near the boundary $\TT\PP^n \setminus \RR^n$ is rather subtle (see \cite{MZ}), but for our purposes (i.e., for faithful tropicalizations), we will not need analysis on the boundary. 

Let $Y^\circ$ be an irreducible variety over $K$, and assume that $Y^\circ$ is embedded in $\GG_m^n$ as a closed subvariety.   The {\em tropicalization} of $Y^\circ$ is 
the subset $\trop\left(Y^{\circ, \an}\right)$ in $\RR^n$. 
By the Bieri--Groves theorem \cite{BG}, 
$\trop\left(Y^{\circ, \an}\right)$ is the support of an integral rational-affine polyhedral 
complex $\Sigma$ in $\RR^n$. 

\subsection{Unimodular tropicalization and faithful tropicalization}
\label{subsec:unimodular:faithful}
Let $X$ be a connected, smooth projective variety over $K$, and let $L$ be an ample line bundle over $X$. Assume that $X$ has a regular strictly semistable model $\Xscr$ of $X$ over $R$. 
Let $S(\Xscr) \subset X^{\an}$ be the skeleton associated to $\Xscr$. 

We have $S(\Xscr) = \bigcup_{S} \Delta_S$, where $S$ runs through all the strata of $\Xscr_s$. Recall that $\Delta_S$ is homeomorphic to $\Delta^{r-1}$ 
for some $r$ by \eqref{eqn:intrinzic}
and that this homeomorphism is intrinsic in the sense that 
it is unique up to reordering of the coordinates. Through this homeomorphism, we will regard $\Delta_S$ as an integral rational-affine polyhedron in $\RR^{r}$. 

Suppose that we are given nonzero global sections $s_0, s_1, \ldots, s_n \in H^0(X, L)$ such that 
\begin{equation}
\label{eqn:map:tropical}
  \varphi: X^\an \longrightarrow \TT\PP^n, \quad p = (p, |\cdot|) \mapsto \left(-\log|s_0(p)|: \cdots: -\log|s_n(p)|\right) 
\end{equation}
is a morphism (i.e., everywhere defined).  Since $X^{\an}$ contains closed points of $X$ as classical points (cf. Example~\ref{eg:classical}), it follows that $\varphi^\prime: X \longrightarrow \PP^{n}$ 
defined by $p  \mapsto \left(s_0(p): \cdots: s_n(p)\right)$ is a morphism. We
denote by the same $\varphi^\prime$ the induced morphism from $X^{\an}$ to $\PP^{n, \an}$. 
Then we have $\varphi = \trop \circ \varphi^\prime$, where 
$\trop:  \PP^{n, \an} \to \TT\PP^n$ is the tropicalization map in \eqref{eqn:trop:map}.

\begin{Remark}
\label{remark:lookinRR}
With the above notation,
we further set  $f_i := s_i / s_0$ for $i = 1 , \ldots , n$.
Then $f_i$ is a nonzero rational function on $X$.
Since any point in $S (\mathscr{X})$ is an
absolute value on the function field of $X$,
we have $|f_i (p)| \neq 0$ for any $p \in S (\mathscr{X})$.
It follows that the restriction of 
$\varphi$ defined in (\ref{eqn:map:tropical}) to $S(\Xscr)$ 
factors through 
\[
  \rest{\varphi}{S(\Xscr)}: S(\Xscr) \to \RR^n, \quad 
  p = (p, |\cdot|) \mapsto (-\log |f_1(p)|, \ldots , -\log |f_n(p)|)  
\]
and the standard embedding
$\RR^n \hookrightarrow \TT\PP^n$, $(x_1, \ldots, x_n) \mapsto (0: x_1: \cdots: x_n)$ in 
\S\ref{subsec:tropical}. 
In the following, we sometimes regard $\rest{\varphi}{S ( \mathscr{X} )}$
as a map valued in $\RR^n$.
\end{Remark}

We put $Y^\circ := \varphi^\prime(X)\cap \GG_m^n \subset \PP^n$. Then $\trop\left(Y^{\circ, \an}\right)$ is the support of some integral rational-affine polyhedral complex $\Sigma$ in $\RR^n$. For $p \in S(\Xscr)$, we have $|s_i(p)| \neq 0$ for any $i$, so that $\varphi^\prime(p) \in Y^{\circ, \an}$. Thus we have $\rest{\varphi}{S(\Xscr)}: S(\Xscr) \to \Sigma \subset \RR^n$. 

\begin{Definition}[cf. {\cite[(5.15)]{BPR}, \cite[Definition~9.2]{GRW}}]
\label{definition:tropicalizations}
Let $s_0, s_1, \ldots, s_n \in H^0(X, L)$ be nonzero global sections of $L$ such that $\varphi: X^\an \longrightarrow \TT\PP^n$ in \eqref{eqn:map:tropical} is a morphism. 
\begin{enumerate}
\item
(unimodular tropicalization) Let $S$ be a stratum of $\Xscr_s$. 
We say that $\varphi$ is {\em unimodular} on $\Delta_S$ if $\Delta_S$ can be covered by finitely many integral rational-affine polyhedra $\Delta$ such that $\rest{\varphi}{\Delta}$ is a unimodular integral rational-affine map on $\Delta$ (cf. \S\ref{subsec:tropical}). 
We say that $\varphi$ is a {\em unimodular tropicalization} of $S(\Xscr)$ if 
$\varphi$ is unimodular on $\Delta_S$ for any stratum $S$ of $\Xscr_s$. 
\item
(faithful tropicalization) We say that $\varphi$ is a {\em faithful tropicalization} of $S(\Xscr)$ if $\varphi$ is a unimodular  tropicalization of $S(\Xscr)$ such that the restriction of 
$\varphi$ to $S(\Xscr)$ is injective. 
\end{enumerate}
\end{Definition}

We say that the skeleton $S(\Xscr)$ {\em has a unimodular \textup{(}resp. faithful\textup{)} tropicalization associated to the linear system $|L|$} if there exist nonzero global sections $s_0, s_1, \ldots, s_n \in H^0(X, L)$ such that the associated morphism $\varphi$ gives a unimodular (resp. faithful) tropicalization of $S(\Xscr)$. 

\begin{Remark}
For the definitions of unimodular and faithful tropicalizations in Definition~\ref{definition:tropicalizations}, we require that $\varphi$ in \eqref{eqn:map:tropical} be everywhere defined. This requirement could be relaxed as long as $\rest{\varphi}{S(\Xscr)}$ is everywhere defined. Suppose that we are given nonzero global sections 
$s_0, s_1, \ldots, s_n \in H^0(X, L)$, and we set 
\begin{equation}
\label{eqn:varphi:dasharrow}
  \varphi: X^\an \dasharrow \TT\PP^n, \quad p = (p, |\cdot|) \mapsto \left(-\log|s_0(p)|: \cdots: -\log|s_n(p)|\right), 
\end{equation}
which may not be defined at some classical points (cf. Example~\ref{eg:classical}). As we argue in Remark~\ref{remark:lookinRR}, $\rest{\varphi}{S(\Xscr)}: S(\Xscr) \to \TT\PP^n$ is everywhere defined. We 
call $\varphi$ in \eqref{eqn:varphi:dasharrow} a {\em rational} unimodular (resp. faithful) tropicalization of $S(\Xscr)$ if $\rest{\varphi}{S(\Xscr)}$ satisfies the condition (1) (resp. (2)) of Definition~\ref{definition:tropicalizations}, and, in this case. we say that 
$S(\Xscr)$ has a {\em rational} unimodular (resp. faithful) tropicalization associated to the linear system $|L|$. 

Since $X^{\an}$ contains closed points of $X$ as classical points, 
we see that $S(\Xscr)$ has a unimodular (resp. faithful) tropicalization associated to 
$|L|$ if and only if  (i) $S(\Xscr)$ has a rational unimodular (resp. faithful) tropicalization associated to $|L|$ and (ii) $L$ is basepoint free. 
\end{Remark}

\section{Vanishing and basepoint-freeness}
\label{sec:lemmas}

In this section, we give some remarks on vanishing of cohomologies
and basepoint-freeness over reducible varieties.

First, we recall a Kodaira-type vanishing theorem for a 
projective 
strictly normal crossing variety. 
Here a {\em strictly} normal crossing variety is 
a normal crossing variety such that each irreducible component is smooth.

\begin{Proposition}
\label{prop:vanishing}
Let $k$ be a field of characteristic $0$.
Let $Z$ be a projective strictly normal crossing variety over $k$.
Let $N$ be an ample line bundle over $X$, and let $\omega_Z$ denote 
the canonical line bundle of $Z$. Then we have 
$H^i(Z, N\otimes \omega_Z) = 0$ for any $i >0$. 
\end{Proposition}

\Proof
By a flat base change, we may assume that $k$ is algebraically closed. 
The assertion follows from, for example, \cite[Theorem~1.1]{Fu} (applied to $Y := \Xscr_s, \Delta := 0$ and $X = \Spec(k)$ in {\it op.\! cit.}). Or, since normal crossing varieties are Cohen-Macaulay, 
the assertion also follows from \cite[Corollary~6.6]{KSS}, where the authors of \cite{KSS} prove the vanishing of the dual $H^i(Z, \omega_Z^{-1}) = 0$ for any $i < \dim Z$.
\QED

Next, we consider basepoint freeness on 
strictly normal crossing varieties.
For a positive integer $d$, we consider 
the following quantity related to the adjoint line bundle over {\em strictly normal crossing varieties} in characteristic zero: 
\begin{equation}
\label{eqn:def:tilde:phi:d}
\widetilde{\phi}(d) 
:= 
\left\{
m_0 \in \ZZ \;\left|\;  
\begin{aligned}
& \text{For any projective strictly normal crossing variety $Z$ with} \\
& \text{$\dim Z \leq d$ over a field of characteristic zero, and for any }\\
& \text{ample line bundle $N$ over $Z$, $N^{\otimes m} \otimes \omega_Z$ is basepoint free}\\
& \text{for any $m \geq m_0$.}
\end{aligned}
\right.
\right\}. 
\end{equation}

Compared with $\phi(d)$ in \eqref{eqn:def:phi:d},  
we allow that $Z$ is 
a projective strictly normal crossing variety;
$Z$ is not necessarily smooth.
By convention,
if there does not exist such an $m_0$, then we set 
$\widetilde{\phi}(d) := + \infty$,
but this does not occur as we see below.

\begin{Lemma}
\label{lem:main:3}
We have $\phi(d) = \widetilde{\phi}(d)$. 
\end{Lemma}

\Proof
It suffices to show that $\widetilde{\phi}(d) \leq \phi(d)$. 
Let $Z$ be a  projective strictly normal crossing variety of dimension less than or equal to $d$ 
 over a field $k$ of characteristic zero, let $N$ be an ample line bundle over $Z$, 
and let $m \geq \phi(d)$. 
We are going to show that $N^{\otimes m} \otimes \omega_Z$ is basepoint free by induction on $\dim Z$. By a flat base change, we may assume that $k$ is algebraically closed. 
Also we may assume that $Z$ is connected. 

We put $\dim Z = e \;(\leq d)$. If $e = 0$, then any line bundle is basepoint free, and there is nothing to prove. Let $e \geq 1$, and let $Z = \bigcup_{i=1}^\ell Z_i$ be the irreducible decomposition of $Z$. 
We take any closed point $p$ of $Z$. Without loss of generality, we may assume that $p$ lies on $Z_1$. We set 
$Z_1^\prime : = \bigcup_{i=2}^\ell Z_i$, and $W_1 := Z_1 \cap Z_1^\prime$.
Note that $W_1$ is a 
projective strictly normal crossing variety of dimension $e-1$. 

\medskip
{\bf Case 1.}\quad  
Suppose that $p$ is a smooth point of $Z_1$. Since $\rest{N^{\otimes m} \otimes \omega_Z}{Z_1}(-W_1) = \rest{N}{Z_1}^{\otimes m} \otimes \omega_{Z_1}$, by the definition of $\phi(d)$, there exists 
$t \in H^0\left(Z_1, \rest{N^{\otimes m} \otimes \omega_Z}{Z_1}(-W_1)\right)$ with $t(p) \neq 0$. Via the natural injection 
$\rest{N^{\otimes m} \otimes \omega_Z}{Z_1}(-W_1) \hookrightarrow \rest{N^{\otimes m} \otimes \omega_Z}{Z_1}$, let $s \in H^0\left(Z_1, \rest{N^{\otimes m} \otimes \omega_Z}{Z_1}\right)$ be the image of $t$. Then $s(p) \neq 0$ and $s$ is equal to zero along $W_1$. Let $\widetilde{s} \in H^0\left(Z, N^{\otimes m} \otimes \omega_Z\right)$ be 
the zero extension of $s$ to $Z$,
that is, $\rest{\widetilde{s}}{Z_1} = s$ and 
$\widetilde{s}$ vanishes outside of $Z_1$. 
Then $\widetilde{s}(p) \neq 0$, which
shows that $N^{\otimes m} \otimes \omega_Z$ is free at $p$. 

\medskip
{\bf Case 2.}\quad  
Suppose that $p$ lies on  $W_1$. By the induction hypothesis, there exists 
$s \in H^0\left(W_1, \rest{N}{W_1}^{\otimes m} \otimes \omega_{W_1} \right)$ such that  
$s(p)\neq 0$.  
By the adjunction formula,
there exists a natural exact sequence
\[
 0 \to \rest{N}{Z_1}^{\otimes m} \otimes \omega_{Z_1} \to \rest{N}{Z_1}^{\otimes m} \otimes \omega_{Z_1}(W_1) \to 
 \rest{N}{W_1}^{\otimes m} \otimes \omega_{W_1} \to 0
,
\]
and by the Kodaira vanishing we have
$h^1(Z_1, \rest{N}{Z_1}^{\otimes m} \otimes \omega_{Z_1} ) = 0$.
It follows that there exists $s_1 \in 
H^0(Z_1, \rest{N}{Z_1}^{\otimes m} \otimes \omega_{Z_1}(W_1))$ with $\rest{s_1}{W_1} = s$.  Similarly, 
we have an exact sequence 
\[
 0 \to \rest{N}{Z_1^\prime}^{\otimes m} \otimes \omega_{Z_1^\prime} \to \rest{N}{Z_1^\prime}^{\otimes m} \otimes \omega_{Z_1^\prime}(W_1) \to 
 \rest{N}{W_1}^{\otimes m} \otimes \omega_{W_1} \to 0
\]
by the adjunction formula and 
the vanishing $h^1(Z_1^\prime, \rest{N}{Z_1^\prime} \otimes \omega_{Z_1^\prime} ) = 0$ by Proposition~\ref{prop:vanishing},
so that there exists $s_1^\prime \in 
H^0(Z_1^\prime, \rest{N}{Z_1^\prime} \otimes \omega_{Z_1^\prime}(W_1))$ with $\rest{s_1^\prime}{W_1} = s$.  Since 
\[
\rest{N}{Z_1}^{\otimes m} \otimes \omega_{Z_1}(W_1) \cong 
\rest{(N^{\otimes m} \otimes \omega_{Z})}{Z_1}
\qquad\text{and}\qquad 
\rest{N}{Z_1^\prime}^{\otimes m} \otimes \omega_{Z_1^\prime}(W_1) \cong 
\rest{(N^{\otimes m} \otimes \omega_{Z})}{Z_1^\prime}, 
\]
we regard $s_1$ and $s_1^\prime$ as sections of $H^0(Z_1, \rest{(N^{\otimes m} \otimes \omega_{Z})}{Z_1})$ and $H^0(Z_1^\prime, \rest{(N^{\otimes m} \otimes \omega_{Z})}{Z_1^\prime})$, respectively. Since $\rest{s_1}{W_1} = \rest{s_1^\prime}{W_1}$, $s_1$ and $s_1^\prime$ glue together to give a section $\widetilde{s} \in H^0(Z, N^{\otimes m} \otimes \omega_{Z})$. 
Then $\widetilde{s}(p) = s(p) \neq 0$, 
which shows that $N^{\otimes m} \otimes \omega_{Z}$ is free at $p$. 
\QED

\section{Unimodular tropicalization}
\label{sec:proof:unimodular}
As before, let $X$ be a connected, smooth projective variety of dimension $d$
over $K$. Let $L$ be an ample line bundle over $X$.
Let $\Xscr$ be a regular strictly semistable model of $X$ over $R$. 
Let $\mathscr{L}$ be a line bundle over $\Xscr$ with 
$\rest{\mathscr{L}}{X} = L$ and assume that $\mathscr{L}$ is relatively ample.

In this section, 
under the assumption in Theorem~\textup{\ref{thm:main}}, 
we show that 
there exist global sections $s_0, s_1, \ldots, s_n$ of 
$L^{\otimes m} \otimes \omega_X$ that give a 
{\em unimodular} tropicalization of the skeleton $S(\Xscr)$.

We prove two technical lemmas.

\begin{Lemma}
\label{lemma:main:2}
Under the assumption of Theorem~\textup{\ref{thm:main}}, 
the restriction map $H^0(\Xscr, \Lscr^{\otimes m}\otimes \omega_{\Xscr/R}) \to 
H^0(\Xscr_s, \Lscr^{\otimes m}\otimes \omega_{\Xscr/R}\vert_{\Xscr_s})$ is surjective. 
\end{Lemma}

\Proof
We note that $\rest{\Lscr^{\otimes m}\otimes \omega_{\Xscr/R}}{\Xscr_s} =  (\rest{\Lscr}{\Xscr_s})^{\otimes m}\otimes \omega_{\Xscr_s}$. 
By the base change theorem
(see \cite[III, 12]{Ha}), it then suffices to show that 
$h^1\left(\Xscr_s, (\rest{\Lscr}{\Xscr_s})^{\otimes m}\otimes \omega_{\Xscr_s}\right) = 0$. 
Since $\Xscr_s$ is a strictly normal crossing variety
and 
$\rest{\Lscr}{\Xscr_s}$ is an ample line bundle over $\Xscr_s$, the vanishing 
$h^1\left(\Xscr_s, (\rest{\Lscr}{\Xscr_s})^{\otimes m}\otimes \omega_{\Xscr_s}\right) = 0$ follows from Proposition~\ref{prop:vanishing}. 
\QED

We fix the notation.
Let  ${\Xscr}_s = \bigcup_{1 \leq i \leq \ell}\Xscr_{s, i}$ denote the irreducible decomposition of 
the special fiber $\Xscr_s$. 
For each $i$ $(1 \leq i \leq \ell)$, let $\bigcup_{j=1}^{b_i} Y_{i j}$ denote the irreducible decomposition of 
$\Xscr_{s, i}\cap (\Xscr_s - \Xscr_{s, i})
=
\Xscr_{s, i} \cap 
\bigcup_{1 \leq i^\prime \leq \ell,\, i^\prime \neq i} \Xscr_{s, i^\prime}$. 

Recall that $\phi (d)$ is a positive integer defined in \eqref{eqn:def:phi:d}. 

\begin{Lemma}
\label{lemma:main:1}
Under the assumption of Theorem~\textup{\ref{thm:main}}, 
$(\Lscr^{\otimes m} \otimes \omega_{\Xscr/R})\vert_{\Xscr_{s, i}}\left(-\sum_{j=1}^{b_i} Y_{i j}\right)$ is basepoint free for any $m \geq \phi(d)$. 
\end{Lemma}

\Proof
By the adjunction formula, we get 
\begin{align*}
(\Lscr^{\otimes m} \otimes \omega_{\Xscr/R})\vert_{\Xscr_{s, i}}\left(-{\textstyle\sum}_{j=1}^{b_i} Y_{i j}\right) 
& \cong 
(\rest{\Lscr}{\Xscr_{s, i}})^{\otimes m} \otimes \rest{\omega_{\Xscr_s}}{\Xscr_{s, i}}\left(-{\textstyle\sum}_{j=1}^{b_i} Y_{i j}\right) \\
& \cong 
(\rest{\Lscr}{\Xscr_{s, i}})^{\otimes m} \otimes \omega_{\Xscr_{s, i}}. 
\end{align*}

Since $\Lscr$ is assumed to be ample, 
$\rest{\Lscr}{\Xscr_{s, i}}$ is ample. 
By the definition of $\phi(d)$, 
$(\rest{\Lscr}{\Xscr_{s, i}})^{\otimes m} \otimes \omega_{\Xscr_{s, i}}$ is basepoint free for 
any $m \geq \phi(d)$. 
Thus 
$(\Lscr^{\otimes m} \otimes \omega_{\Xscr/R})\vert_{\Xscr_{s, i}}\left(-{\textstyle\sum}_{j=1}^{b_i} Y_{i j}\right)$ is basepoint free for 
any $m \geq \phi(d)$. 
\QED

\subsection*{Well-behaved global sections}
By Lemma~\ref{lemma:main:1}, there exists a non-zero global section 
\[
\xi_i \in H^0\left(
\Xscr_{s, i}, \rest{\left(\Lscr^{\otimes m} \otimes \omega_{\Xscr/R}\right)}{\Xscr_{s, i}}\left(-{\textstyle\sum}_{j=1}^{b_{i}} Y_{ij}\right)
\right)
\]
that does not vanish at the generic point of 
any minimal stratum of $\Xscr_{s, i}$. 
We denote by 
$
  \eta_i^\prime \in H^0\left(
\Xscr_{s, i}, \rest{\left(\Lscr^{\otimes m} \otimes \omega_{\Xscr/R}\right)}{\Xscr_{s, i}}\right)
$
the image of $\xi_i$ under the natural inclusion 
\[
\rest{\left(\Lscr^{\otimes m} \otimes \omega_{\Xscr/R}\right)}{\Xscr_{s, i}}(-{\textstyle\sum}_{j=1}^{b_{i}} Y_{ij}) \hookrightarrow 
\rest{\left(\Lscr^{\otimes m} \otimes \omega_{\Xscr/R}\right)}{\Xscr_{s, i}}. 
\]
We then have $\ord_{Y_{ij}}(\eta_i^\prime) = 1$ for any $j = 1, \ldots, b_i$. 
It follows that there exists a
global section $\eta_i$ of
$\rest{\left(\Lscr^{\otimes m} \otimes \omega_{\Xscr/R}\right)}{\Xscr_s}$
that is the zero-extension of $\eta_i^\prime$,
that is, 
$\rest{\eta_i}{\Xscr_{s, i}} = \eta_i^\prime$
and $\rest{\eta_i}{\Xscr_{s, i'}} = 0$ for $i' \neq i$.
\if
the zero extension of $\eta_i^\prime$ to $\Xscr_s$. Thus we obtain a 
global section 
\[
 \eta_i \in H^0\left(
\Xscr_s, \rest{\left(\Lscr^{\otimes m} \otimes \omega_{\Xscr/R}\right)}{\Xscr_s}\right). 
\]
\fi
By Lemma~\ref{lemma:main:2}, there exists $\widetilde{\eta}_i \in H^0(\Xscr, \Lscr^{\otimes m} \otimes \omega_{\Xscr/R})$ such that $\rest{\widetilde{\eta}_i}{\Xscr_s} = \eta_i$. 

Let $\Xscr_{s, i^\prime} \in \Irr(\Xscr_s)$ ($i^\prime \neq i$) be an irreducible component of $\Xscr_s$ such that 
$\Xscr_{s, i} \cap \Xscr_{s, i^\prime} \neq \emptyset$.  Since $\eta_i$ is the zero extension of 
$\eta_i^\prime$, we have $\ord_{\Xscr_{s, i^\prime}}(\widetilde{\eta}_{i}) \geq 1$. 
We claim that 
\begin{equation}
\label{eqn:order}
\ord_{\Xscr_{s, i^\prime}}(\widetilde{\eta}_{i}) = 1.
\end{equation} 
Indeed, let $Y_{ii'}$ be an irreducible component of $\Xscr_{s, i} \cap \Xscr_{s, i^\prime}$. 
Let $\Vscr$ be a small Zariski open neighborhood of the generic point of $Y_{ii'}$  in $\Xscr$. 
Let $T_{i^\prime} = 0$ be a local equation of $\Xscr_{s, i^\prime}$ on $\Vscr$. 
We write 
$\rest{\widetilde{\eta}_{i}}{\Vscr} = T_{i^\prime}^a f_{i^\prime}$ with $a \geq 0$ and 
$\rest{f_{i^\prime}}{\Xscr_{s, i^\prime}} \not\equiv 0$ on $\Vscr$.
Since $\rest{\widetilde{\eta}_{i}}{\Xscr_{s, i^\prime}} = 0$,
we have $a \geq 1$.
On the other hand,
since $\rest{\widetilde{\eta}_{i}}{\Xscr_{s, i}} = 
\rest{T_{i^\prime}}{\Xscr_{s, i}}^a \rest{ f_{i^\prime}}{\Xscr_{s, i}}$, we have 
$\ord_{Y_{ii'}}(\rest{\widetilde{\eta_i}}{\Xscr_{s, i}}) \geq a$. Since $\ord_{Y_{ii'}}(\rest{\widetilde{\eta_i}}{\Xscr_{s, i}}) = \ord_{Y_{ii'}}(\eta_i^\prime) = 1$, we obtain $1 \geq a$. 
Thus $a = 1$, and we obtain \eqref{eqn:order}. 

We set 
\begin{equation}
\label{eqn:s:i}
  s_i := \rest{\widetilde{\eta}_i}{X} \in H^0\left(X, L^{\otimes m}\otimes \omega_X\right). 
\end{equation}

\subsection*{Base global section}
Since $m \geq \phi(d)$ and since 
$\phi(d) = \widetilde{\phi}(d)$ by Lemma~\ref{lem:main:3}, $(\rest{\Lscr}{\Xscr_s})^{\otimes m}\otimes \omega_{\Xscr_s} = \rest{\left(\Lscr^{\otimes m} \otimes \omega_{\Xscr/R}\right)}{\Xscr_s}$ is basepoint free. 
It follows that there exists a global section 
\[
\eta_0 \in H^0\left(
\Xscr_s, \rest{\left(\Lscr^{\otimes m} \otimes \omega_{\Xscr/R}\right)}{\Xscr_s}\right)
\]
such that $\eta_0$ does not vanish at the generic point of any minimal stratum of $\Xscr_s$.  
By Lemma~\ref{lemma:main:2}, there exists $\widetilde{\eta}_0 \in H^0(\Xscr, \Lscr^{\otimes m} \otimes \omega_{\Xscr/R})$ such that $\rest{\widetilde{\eta}_0}{\Xscr_s} = \eta_0$. 
We set 
\begin{equation}
\label{eqn:s:0}
  s_0 := \rest{\widetilde{\eta}_0}{X} \in H^0\left(X, L^{\otimes m}\otimes \omega_X\right). 
\end{equation}

\subsection*{Unimodular tropicalization}
Let $s_0, s_1, \ldots, s_\ell$ be global sections of $L^{\otimes m}\otimes \omega_X$ constructed in \eqref{eqn:s:i} and \eqref{eqn:s:0}. Since $m \geq \phi(d)$, the definition of $\phi(d)$ tells us that there exist global sections 
$s_{\ell+1}, \ldots, s_{\ell+d + 1}$ of $L^{\otimes m}\otimes \omega_X$ with 
$\bigcap_{i=1}^{\ell+d+1} \zero(s_{i}) = \emptyset$. 
Then $X \to \PP^{\ell+d+1}$ given by $p \mapsto (s_0(p): \cdots : s_{\ell+d}(p))$ is a morphism. 
Let 
\begin{equation}
\label{eqn:unimodular:initial}
\varphi: X^{\an} \to \TT\PP^{\ell+d+1}, \quad
p = (p, |\cdot|) \mapsto (-\log|s_0(p)|: \cdots : -\log|s_{\ell+d+1}(p)|)
\end{equation}
be the associated morphism. 
For $i = 1, \ldots, \ell$, 
we define the rational function $f_i \in \Rat(X) = \Rat(\Xscr)$ by 
\[
  f_i := \frac{s_i}{s_0} = \rest{\frac{\widetilde{\eta}_i}{\widetilde{\eta}_0}}{X}. 
\]
As noted in Remark~\ref{remark:lookinRR},
we identify $\rest{\varphi}{S(\Xscr)}$
with
\begin{align}
\label{eqn:unimodular:initial2}
\rest{\varphi}{S(\Xscr)}: S(\Xscr) \to \RR^{\ell + d + 1}, \quad 
p \mapsto
\left(
- \log |f_1 (p)|
,
\ldots 
,
- \log |f_{\ell + d + 1 } (p)|
\right)
.
\end{align}

Let us
prove that the map in (\ref{eqn:unimodular:initial2}) 
is a unimodular tropicalization
for $S(\Xscr)$. 
We take any minimal stratum $S$ of $\Xscr_s$.
Let $\Delta_S$ be the canonical simplex 
associated to $S$ in $S(\Xscr)$. 
By the construction of the skeleton $S(\Xscr)$, it suffices to check that 
$\rest{\varphi}{\Delta_S}$ is a unimodular integral rational-affine map 
on $\Delta_S$. 
Let $v_{j_1}, \ldots, v_{j_r}$ be the vertices of $\Delta_S$, 
where $v_{j}$ is the Shilov point of $\Xscr_{s, j}$. 
We set $J := \{j_1, \ldots, j_r\} \subset \{1, \ldots, \ell\}$. 
We take a sufficiently small Zariski open neighborhood 
$\Uscr$
of the generic 
point of $S$ in $\Xscr$. 
Then,
since $\widetilde{\eta}_0$ does not vanish at any minimal stratum
and
hence does not at any stratum, we have 
$\rest{\zero(\widetilde{\eta}_0)}{\Uscr} = 0$. 
Further, for $j \in J$, 
by \eqref{eqn:order}, we have 
\[
\rest{\zero(\widetilde{\eta}_j)}{\Uscr} = 
\sum_{
 \substack{
j^\prime \in J\\ j^\prime \neq j}} \Xscr_{s, j^\prime}. 
\]
It follows that,
regarding $f_j$ as an element of $\Rat(\Xscr)$, we have 
\begin{equation}
\label{eqn:fi:order}
\rest{\zero(f_j)}{\Uscr} = \sum_{ \substack{
j^\prime \in J\\ j^\prime \neq j}} \Xscr_{s, j^\prime}. 
\end{equation}

We consider the morphism
\begin{align}
\label{eqn:unimodular2}
  \psi_S: \Delta_S \to \RR^{r}, \quad 
  p = (p, |\cdot|) \mapsto  (-\log|f_{j_1}(p)|, \ldots, -\log|f_{j_r}(p)|).  
\end{align}
We study the function $-\log|f_{j_a}|$ for $a = 1, \ldots, r$ in more detail. 
For $j_a \in J$, let $T_{j_a} = 0$ be a local equation of $\Xscr_{s, j_{a}}$ at the generic point of $S$.  
By \eqref{eqn:fi:order}, we write 
$f_{j_a} = \lambda \prod_{1 \leq b \leq r, b \neq a} T_{j_b}$ for some invertible function $\lambda$ on $\Uscr$. 
We take any $p = (\eta, |\cdot|_{\mathbf{u}, S}) \in \Delta_S$,
with the notation in 
\eqref{eqn:coord:DS}. 
We write 
$\mathbf{u} = (u_1 , \ldots , u_r)$. 
Then by the definition of $|\cdot|_{\mathbf{u}, S}$, we get 
\[
-\log|f_{j_a}(p)|  
= \sum_{1 \leq b \leq r, b \neq a} u_b. 
\]
Since $p$ corresponds to $\mathbf{u} = (u_1 , \ldots , u_r)$ 
via the isomorphism $\Delta_S \cong \Delta^{r-1}$
of simplices, we see 
in particular that $-\log|f_{j_a}(p)|$ is an affine function on $\Delta_S$. 
Thus $\psi_S$ is an affine map.

Put $\mathbf{v}_a = \psi_S(v_{j_a})$ for $a = 1, \ldots , r$. Then $\psi_S(\Delta_S)$ is the $(r-1)$-dimensional simplex spanned by $\mathbf{v}_1, \mathbf{v}_2, \ldots, \mathbf{v}_{r}$. 
To show that $\psi_S: \Delta_S \to \RR^r$ is a unimodular integral rational-affine map,  
it suffices to verify that $\mathbf{v}_2 - \mathbf{v}_1, \ldots, \mathbf{v}_{r} - \mathbf{v}_1$ is a part of 
a $\ZZ$-basis of $\ZZ^{r} \subset \RR^r$.  
Put $\mathbf{1} = (1, 1, \ldots, 1)$, and $\mathbf{e}_1 = (1, 0, \ldots, 0), 
\ldots, \mathbf{e}_{r} = (0, 0, \ldots, 1)$ in $\ZZ^{r}$.  
Since
$\mathbf{v}_a = \mathbf{1} - \mathbf{e}_a$ for $a = 1, \ldots, r$, 
$\mathbf{v}_a - \mathbf{v}_1 = \mathbf{e}_1 - \mathbf{e}_a$ for  $a = 2, \ldots, r$. 
Since $\mathbf{e}_1, 
\mathbf{e}_1 - \mathbf{e}_2 ,
\ldots, 
\mathbf{e}_1 - \mathbf{e}_r
$ 
is a $\ZZ$-basis of $\ZZ^{r}$, 
it follows that $\mathbf{v}_2 - \mathbf{v}_1, \ldots, \mathbf{v}_{r} - \mathbf{v}_1$ is a part of a $\ZZ$-basis of $\ZZ^{r}$.  
Thus we have shown that $\psi_S$ in \eqref{eqn:unimodular2} is unimodular on $\Delta_S$. Since
$\psi_S$ factors through
$\rest{\varphi}{\Delta_S}$ , transitivity of lattice indices show that 
$\varphi$ in \eqref{eqn:unimodular:initial2} is unimodular on $\Delta_S$ (see \cite[Lemma~5.17]{BPR} and \cite[Lemma~9.3]{GRW}). 
This concludes that $\varphi$ is a unimodular tropicalization of $S(\Xscr)$. 

\begin{Remark}
\label{remark:inj:on:simplex}
In the above, we have shown that $\psi_S$ in \eqref{eqn:unimodular2}
is affine and unimodular.
It follows that $\psi_S$ is injective.
Since $\psi_S$ factors through
$\rest{\varphi}{\Delta_S}$,
this shows that $\rest{\varphi}{\Delta_S}$ is injective.
\end{Remark}

We put together some properties of $-\log|f_{j_a}|$ on $S(\Xscr)$, which we use in the next section. 

\begin{Lemma}
\label{lemma:f:i}
As above,
let $v_{j_1}, \ldots, v_{j_r}$ be the vertices of $\Delta_S$. For $a = 1, \ldots, r$, 
we set $g := -\log|f_{j_a}|$. 
Then $g$ is affine on $\Delta_S$. 
Further, for any vertex $v$ of $S(\Xscr)$, we have 
\[
g(v) 
\begin{cases}
= 0 & (\text{if $v = v_{j_a}$}),\\
= 1 & (\text{if $v = v_{j_b}$ for some $b = 1, \ldots, r$ with $b \neq a$}).\\
\geq 1 & (\text{otherwise}). 
\end{cases}
\]
\end{Lemma}

\Proof
We have already seen that $g$ is affine on $\Delta_S$.  
Recall that $f_{j_a} = \widetilde{\eta}_{j_a}/\widetilde{\eta}_0$ as a rational function on $\Xscr$ and that $\widetilde{\eta}_0$ does not vanish at the generic point of any stratum. 
Recall also that $v_{j_a}$ is the Shilov point of $\mathscr{X}_{s , j_a}$.
Since $\widetilde{\eta}_{j_a}$ does not vanish at the generic point of $\Xscr_{s, j_a}$, the first equality follows. 
The second equality 
follows from 
\eqref{eqn:order}. 
Further, the third inequality follows from 
the fact that $\rest{\widetilde{\eta}_{j_a}}{\mathscr{X}_{s , i'}} = 0$
for $i' \neq i$.
\QED

\section{Faithful tropicalization}
\label{sec:proof:faithful}
In this section, we show that the unimodular tropicalization of $S(\Xscr)$ constructed 
in \eqref{eqn:unimodular:initial} is actually a faithful tropicalization.
We only need to show that $\varphi$ is {\em injective}. 

We keep the notation in \S\ref{sec:proof:unimodular}. 
We begin with some lemmas. 

\begin{Lemma} 
\label{lem:concavity}
Let $S$ be a stratum of $\mathscr{X}_s$.
Let $\phi$ be a non-zero rational function on $\mathscr{X}$.
Let $\Delta_S \subset X^{\an}$ be the canonical simplex associated to $S$, 
and we set $g := -\log | \phi | : \Delta_S \to \RR$.
Suppose that 
there exist a Zariski open neighborhood 
$\Uscr$ in $\Xscr$
of the generic point of $S$ 
and
a vertical divisor $V$ on $\mathscr{U}$
such that
of $\zero \left( \rest{\phi}{\mathscr{U}} \right) - V$ 
is a horizontal effective divisor on $\Uscr$. 
Let $v_{j_1} , \ldots , v_{j_r}$ be the vertices of $\Delta_S$, where $v_j$ is the Shilov point of 
$\Xscr_{s, j}$. 
For $\mathbf{u} = (u_1 , \ldots , u_r) \in \Delta^{r-1}$, we take the point  
$p = (\eta, |\cdot|_{\mathbf{u}, S}) \in \Delta_S$.
Then, we have
\begin{equation}
\label{eqn:g:ineq}
g (p)
\geq
u_1 g (v_{j_1}) + \cdots + u_r g (v_{j_r})
.
\end{equation}
\end{Lemma}

\Proof
As before, let $T_{j} = 0$ be a local equation of $\Xscr_{s, j}$ on $\Uscr$. 
Shrinking $\Uscr$ if necessary,  
we write $V = \sum_{a=1}^r m_a \Xscr_{s, j_a}$, and  
set $\phi^\prime = T_{j_1}^{-m_1} \cdots T_{j_r}^{-m_r}$, so that 
$\zero ( \phi' ) = - V$ on $\Uscr$. 
Then on $\Uscr$ we have
\[
- \log | \phi^\prime |
= m_1( \log |T_{j_1}|) + \cdots + 
m_r ( \log | T_{j_r}|). 
\]
Since $ \log |T_{j_1}|, \ldots ,  \log |T_{j_r}|$
are affine functions on $\Delta_S$, we see that 
$h:= - \log | \phi^\prime |$ is an affine function on $\Delta_S$. 

By assumption, 
$\zero ( \phi \phi') = \zero(\phi) - V$ is a horizontal effective divisor on $\Uscr$.
Set $g' := - \log | \phi \phi '|$.
Then
$
g' = g+ h
$.
Since $h$ is affine,
\[
h (p)
=
u_1 
h (v_{j_1}) + \cdots + 
u_r 
h (v_{j_r})
.
\]
This means that
if we show the inequality for $g'$ in place of $g$ in \eqref{eqn:g:ineq}, 
then the same inequality holds also for $g$.
Thus we may and do 
assume that $\zero (\phi)$ is a horizontal effective divisor on $\Uscr$.

Since $\ord_{{\Xscr_{s, j_a}}}(\phi) = 0$, we get 
$g (v_{j_a}) = 0$ for any $a = 1, \ldots, r$.
Thus it suffices to show that 
$g (p) \geq 0$.
Let $\xi$ denote the generic point of $S$. 
From the construction of $|\cdot|_{\mathbf{u}, S}$ in 
\eqref{def:mathu:S}, if $f \in \OO_{\Xscr, \xi}$, then 
$|f|_{\mathbf{u}, S} \leq 1$.  Since  $\phi$ is a regular function on $\Uscr$ 
and $\Uscr$ contains $\xi$, we see that $\phi \in  \OO_{\Xscr, \xi}$. 
We then have $|\phi|_{\mathbf{u}, S} \leq 1$. 
This proves $g (p) \geq 0$.
\QED

\begin{Theorem}
The map $\varphi$ in \eqref{eqn:unimodular:initial}
is a faithful tropicalization of $S(\Xscr)$. 
\end{Theorem}

\Proof
Since we have shown that $\varphi$ is a unimodular tropicalization of $S(\Xscr)$ in the previous section, 
we have only to show that $\varphi$ is injective. 

We have $S ( \mathscr{X} ) = \coprod_{S} \relint ( \Delta_S )$, where $S$ runs through 
all the strata of $\Xscr_s$. 
As noted in Remark~\ref{remark:inj:on:simplex}, 
$\varphi$ is injective on each canonical simplex. 
Thus it suffices to show that 
$\varphi ( \relint ( \Delta_S ) ) \cap \varphi ( \relint ( \Delta_T ) )
= \emptyset$ for any strata $S$ and $ T$ with $S \neq T$.
Furthermore,
since it is injective on each canonical simplex,
we may and do assume that $\Delta_T \nsubseteq \Delta_S$ 
and $\Delta_S \nsubseteq \Delta_T$. 
By symmetry, we assume that $\Delta_T \nsubseteq \Delta_S$. 
As before, let $v_{j_1} , \ldots , v_{j_r}$ be the vertices of $\Delta_S$, where $v_j$ is the Shilov point of 
$\Xscr_{s, j}$. 
Since $\Delta_T \nsubseteq \Delta_S$, there exists $1 \leq a \leq r$ such that $v_{j_a} \not\in \Delta_T$.  

We put $f_{j} := s_{j} / s_0$ for all $j$.
Recall from 
Remark~\ref{remark:lookinRR}
that $\varphi$ is identified with the map
\[
S(\Xscr) \to \RR^{\ell + d + 1}, \quad p = (p, |\cdot|) \mapsto 
( -\log|f_{1}(p)| , \ldots , - \log |f_{\ell + d +1} (p)| ). 
\]
We consider the function
\[
g : S(\Xscr) \to \RR , \quad p = (p, |\cdot|) \mapsto -\log|f_{j_a}(p)|. 
\]
It suffices to prove that
$
g( \relint ( \Delta_S ) ) \cap 
g( \relint ( \Delta_T ) )
= \emptyset$. 
First, we take any $p \in \relint ( \Delta_S )$,
and we write
$p = (\eta, |\cdot|_{\mathbf{u}, S})$ for some
$\mathbf{u} = (u_1 , \ldots , u_r) \in \relint(\Delta^{r-1})$
(cf. \eqref{eqn:coord:DS}). 
By Lemma~\ref{lemma:f:i}, we get 
\[
 g(p)
 = u_1 g (v_{j_1}) + 
\cdots + 
u_r g (v_{j_r}) \\
 = \sum_{1\leq b\leq r, b \neq a} u_b.  
\]
Since 
$\mathbf{u} \in \relint(\Delta^{r-1})$, 
we have $u_1 > 0, \ldots, u_r >0$.
Since 
$u_1 + \cdots + u_r = 1$,
it follows that $0< \sum_{1\leq b\leq r, b \neq a} u_b <1$.
Thus $0 < g(p) < 1$.

Next, we take any $q \in \Delta_T$. Let $v_{\ell_1}, \ldots, v_{\ell_t}$ be the vertices of $\Delta_T$. 
By the definition of $\Delta_T$
(cf. \eqref{eqn:coord:DS}), there exist
$\mathbf{u}^\prime = (u_1^\prime, \ldots, u_t^\prime)$ with $u_1^\prime \geq 0, \ldots, u_t^\prime \geq 0$ and 
$u_1^\prime + \cdots + u_t^\prime = 1$ such that $q =  (\eta, |\cdot|_{\mathbf{u}^\prime, T})$.
Let $\Uscr^\prime$ be a sufficiently small Zariski open neighborhood of the generic point of $T$. 
Recall that $f_{j_a} = \widetilde{\eta}_{j_a}/\widetilde{\eta}_0$ as a rational function on $\Xscr$ and that $\widetilde{\eta}_0$ does not vanish at the generic point of any stratum. Thus $\rest{\zero(\widetilde{\eta}_0)}{\Uscr^\prime} = 0$. It follows that the horizontal part of $\rest{\zero(f_{j_a})}{\Uscr^\prime}$ is effective. 
We apply Lemma~\ref{lem:concavity}, and then we get 
\begin{align*}
  g(q) \geq u_1^\prime g(v_{\ell_1})  + \cdots + u_t^\prime g(v_{\ell_t}) . 
\end{align*}
Since $v_{\ell_1} \neq v_{j_a} , \ldots , v_{\ell_t} \neq v_{j_a}$,
Lemma~\ref{lemma:f:i} gives us
$g(v_{\ell_1}) \geq 1, \ldots, g(v_{\ell_t}) \geq 1$.
It follows that 
$u_1^\prime g(v_{\ell_1})  + \cdots + u_t^\prime g(v_{\ell_t}) \geq 1$.
Thus we have
$g(q) \geq 1$.

Therefore, we have $g( \relint ( \Delta_S ) ) \subseteq (0,1)$,
and $g(\Delta_T) \subset [1 , + \infty )$.
Thus, in particular,
$g( \relint ( \Delta_S ) ) \cap g(\relint(\Delta_T))
= \emptyset$. 
This completes the proof. 
\QED

\end{document}